\documentclass{article}
\usepackage{amssymb,amsmath,amsthm,graphics}
\newtheorem{theorem}{Theorem}[section]
\newtheorem{lemma}[theorem]{Lemma}
\newtheorem{question}[theorem]{Question}
\newtheorem{conjecture}[theorem]{Conjecture}

\newcommand{\bq}{\[}
\newcommand{\eq}{\]}
\newcommand{\C}{\mathbb{C}}
\newcommand{\dual}{{^\ast}}
\newcommand{\odual}{{^\circ}}
\newcommand{\tensor}{\otimes}
\newcommand{\Inv}{\mathrm{Inv}}
\newcommand{\Hom}{\mathrm{Hom}}
\newcommand{\Span}{\mathrm{Span}}
\newcommand{\bd}{\partial}
\renewcommand{\sp}{\mathrm{sp}}
\newcommand{\osp}{\mathrm{sp}}
\newcommand{\so}{\mathrm{so}}
\renewcommand{\sl}{\mathrm{sl}}

\begin{document}

\title{The Quantum $G_2$ Link Invariant}
\author{Greg Kuperberg \thanks{Supported by a National Science Foundation
graduate fellowship in mathematics and a Sloan Foundation graduate fellowship
in mathematics.}
\\ University of Chicago}
\date{October 7, 1991}

\maketitle

\begin{abstract}
We derive an inductive, combinatorial definition of a polynomial-valued
regular isotopy invariant of links and tangled graphs.  We show that
the invariant equals the Reshetikhin-Turaev invariant corresponding to
the exceptional simple Lie algebra $G_2$.  It is therefore related to
$G_2$ in the same way that the HOMFLY polynomial is related to $A_n$
and the Kauffman polynomial is related to $B_n$, $C_n$, and $D_n$.  We
give parallel constructions for the other rank 2 Lie algebras and
present some combinatorial conjectures motivated by the new inductive
definitions.
\end{abstract}

This paper is divided into two parts.  In the first part we derive from
first principles some variants of the link invariant known as the Jones
polynomial.  In the second part we show that these invariants are
the same known invariants constructed using rank 2 Lie algebras, and
we discuss some of their properties.

\section{Invariants of links and graphs}

The simplest known definition of the Jones polynomial is the Kauffman
bracket \cite{Kauffman:state}, which in this paper will be denoted by
$\langle\cdot\rangle_{A_1}$ and will be called the $A_1$ bracket.  The
$A_1$ bracket is given by the following recursive rules: $$\includegraphics{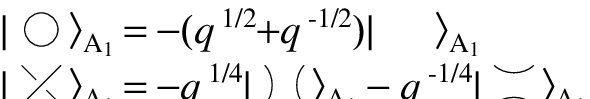}$$

The goal of this part of the paper is to derive definitions of
the following three variants of the Jones polynomial:

\begin{theorem}  There is an invariant for regular isotopy of
projections of links and tangled trivalent graphs called
$\langle\cdot\rangle_{G_2}$ which is given by the following recursive rules:
$$\includegraphics{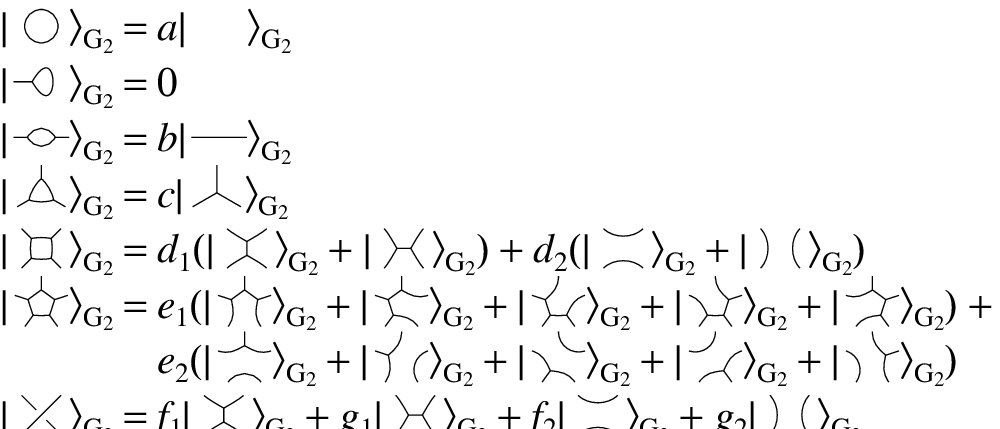}$$
where the coefficients are defined as:
\begin{eqnarray*}
a   & = & q^5 + q^4 + q + 1 + q^{-1} + q^{-4} + q^{-5} \\
b   & = & - q^3 - q^2 - q - q^{-1} - q^{-2} - q^{-3} \\
c   & = & q^2 + 1 + q^{-2} \\
d_1 & = & - q  - q^{-1} \\
d_2 & = & q + 1 + q^{-1} \\
e_1 & = & 1 \\
e_2 & = & -1 \\
f_1 & = & \frac{1}{1 + q^{-1}} \\
g_1 & = & \frac{1}{1 + q} \\
f_2 & = & \frac{q}{1 + q^{-1}} \\
g_2 & = & \frac{q^{-1}}{1 + q},
\end{eqnarray*}
$q$ is an indeterminate, and $\langle\emptyset\rangle_{G_2} = 1$.
\label{existg2} \end{theorem}

\begin{theorem}  There is an invariant for regular isotopy of
projections of links and tangled trivalent graphs called
$\langle\cdot\rangle_{A_2}$ which is given by the following recursive rules:
$$\includegraphics{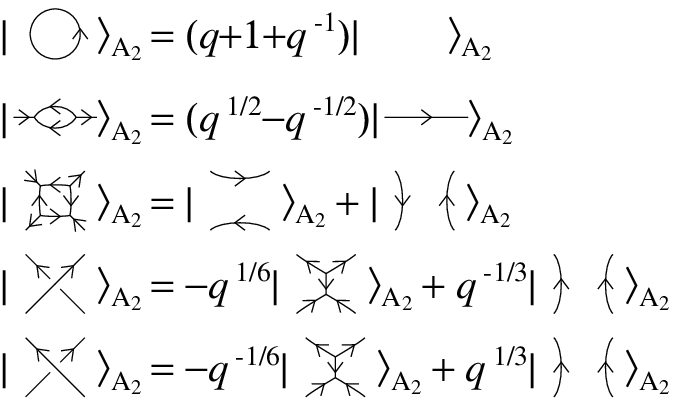}$$
where $q$ is an indeterminate and $\langle\emptyset\rangle_{A_2} = 1$.
\label{exista2} \end{theorem}

\begin{theorem}  There is an invariant for regular isotopy of
projections of links and tangled trivalent graphs called
$\langle\cdot\rangle_{C_2}$ which is given by the following recursive rules:
$$\includegraphics{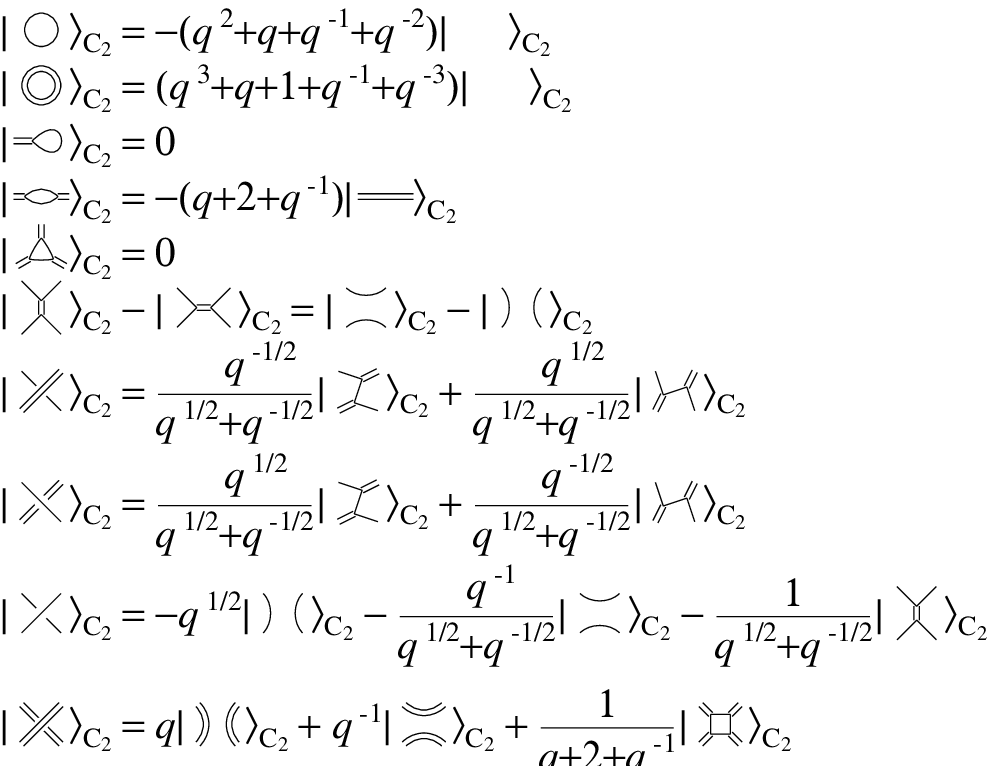}$$
where $q$ is an indeterminate and $\langle\emptyset\rangle_{C_2} = 1$.
Part of the invariant can also be defined by the recursive rules:
$$\includegraphics{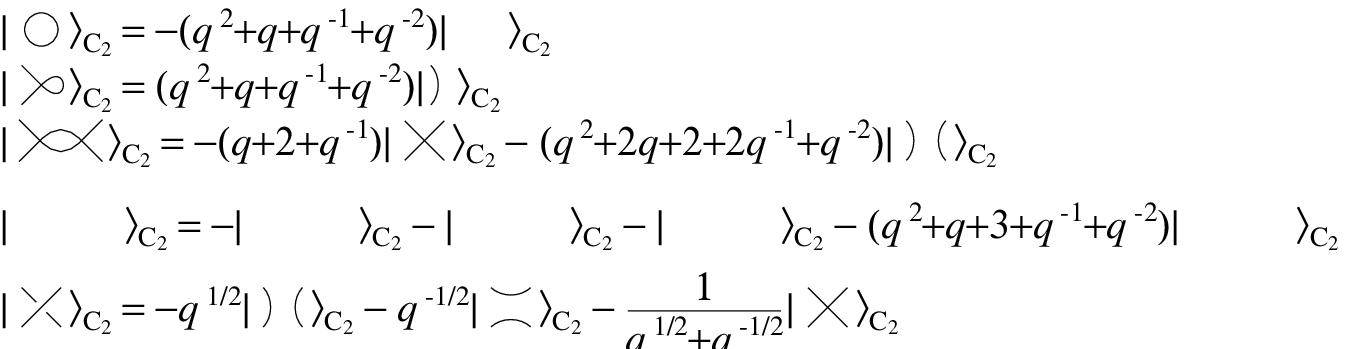}$$
which are related to the first set of rules by: $$\includegraphics{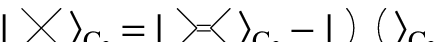}$$
\label{existc2}
\end{theorem}

These invariants are not new.  In fact, all three invariants are
special cases of the Reshetikhin-Turaev invariants
\cite{RT:ribbon}.  (It should be emphasized that these
invariants are part of a wide literature of closely related ideas.  Two
other influential papers are \cite{Drinfeld:quantum} and \cite{Witten:jones}.)
According to Reshetikhin and Turaev, for each simple Lie algebra
$\mathfrak{g}$, there exists an invariant $RT_{\mathfrak{g}}$ of appropriately
colored tangled ribbon graphs.  Each edge is colored by an irreducible
representation of $\mathfrak{g}$ and each vertex should be colored by a
tensor of a certain kind.  The invariants we define all satisfy
$\langle G \rangle_{\mathfrak{g}} = RT_{\mathfrak{g}}(G)$ if the graph is
colored in a certain simple way.  (By our notation, the Kauffman
bracket also satisfies this equation.) This connection will be
explained in detail in the second part of the paper.

In addition, the $A_2$ bracket is essentially a specialization of the
HOMFLY polynomial \cite{HOMFLY:invariant}, because the HOMFLY polynomial describes
$RT_{A_n}$; and the $C_2$ bracket is essentially a specialization of
the Kauffman polynomial, because the Kauffman polynomial 
\cite{Kauffman:regular}

describes $RT_{B_n}$, $RT_{C_n}$, and $RT_{D_n}$.  By the
relation $B_2 = C_2$, the $C_2$ bracket encompasses a second
specialization of the Kauffman polynomial.  (The Kauffman polynomial
should not be confused with the Kauffman bracket.) Thus, the most
interesting case is the $G_2$ bracket.  The Reshetikhin-Turaev
definition in the $G_2$ case is also cast in an explicit form in
\cite{Reshetikhin:g2}.

The author recently learned that many of the results presented here, in
particular the definition of the $A_2$ and $G_2$ brackets, were
obtained independently by Francois Jaeger 
\cite{Jaeger:confluent,Jaeger:bipartite,Jaeger:kauffman}.

\section{Acknowledgments}

The author would like to thank Andrew Casson, Vaughan Jones, Robion Kirby,
Kenneth Millett, Nicolai Reshetikhin, and Oleg Viro for their sustained
interest in the work that resulted in this paper and for their helpful
comments.  The author is also indebted to Rena Zieve for checking certain
important calculations.

\subsection{Topological preliminaries}

By a planar graph we mean a finite combinatorial graph
which is embedded in the sphere $S^2$.  We allow planar graphs to
have edges which go in a circle and have no vertices, as well
as a multiple edges between vertices and vertices connected
to themselves.  A face of
a planar graph is a connected component of the complement of the
graph.  When interpreting the illustrations of planar graphs in this
paper, the reader should identify the border of the page to a point.

Let $D$ be a disk in $S^2$.  We define a planar graph with boundary
to be the intersection of $D$ with a planar graph which is transverse
to the boundary of $D$.  Thus a planar graph with boundary is a graph
which is embedded in a disk such that the edges of the graph meet the
boundary of the disk transversely at special vertices of degree 1,
called endpoints.  The boundary of a planar graph in $D$ is defined to
be the boundary of $D$ together with the endpoints and an
inward-pointing normal.  If $b$ is a boundary, the opposite boundary
$-b$ is defined to be the same object with the distinguished normal
vector reversed.  In this paper we will consider graphs decorated in
various ways and it is understood that the endpoints of a graph with
boundary should be decorated correspondingly.  For example, if the
edges of a graph are colored, the endpoints should be assigned the same
colors.  If the edges are oriented, each endpoint should be assign a
distinguished normal vector.  The idea is that if two possibly
decorated planar graphs have opposite boundaries, their union should be
a planar graph without boundary decorated in the same way.

We define a graph projection to be a planar graph with
special tetravalent vertices, called crossings, such that one pair
of opposite incoming edges is labeled as passing ``over'' the other
pair, like so: $$\includegraphics{fig7.eps}$$
We define a graph projection with boundary similarly.

We define a regular isotopy of a graph projection to be a sequence
of operations consisting of ambient isotopy of the 2-sphere and the following
combinatorial moves: $$\includegraphics{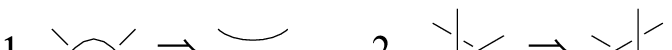}$$
There are versions of the last two moves for every possible kind of
vertex in the graph, although the moves are only illustrated for a
trivalent vertex.  We define a tangled ribbon graph (introduced in
\cite{RT:ribbon}) to be an equivalence class of graphs
projections under regular isotopy.  Regular isotopy moves for graph
projections with boundary are the same, except that an ambient isotopy
must leave the boundary of the graph fixed.  We can also consider
regular isotopy of graphs decorated in various ways; there is a version
of each regular isotopy move for every possible decoration of the arcs
involved in the move.

Two important classes of graph projections are those with no vertices,
which are called link projections, and those with only trivalent
vertices, which we will call freeway projections.  (A link projection
with boundary is allowed to have univalent vertices at the boundary.)
Here is an example of a freeway projection: $$\includegraphics{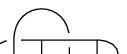}$$
A ribbon link (also called a framed link) is an equivalence class of
link projections, a freeway is an equivalence class of freeway
projections, a ribbon link with boundary (also called a framed tangle)
is an equivalence class of ribbon links with boundary, and a freeway
with boundary is an equivalence class of freeway projections with
boundary.  Here is an example of a freeway projection with boundary:
$$\includegraphics{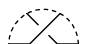}$$

A tangled ribbon graph has a natural 3-dimensional interpretation as
(the equivalence class under ambient isotopy of) a graph embedded in
$S^3$ together with a (continuous) distinguished normal vector at each
point.  Alternatively we can think of a tangled ribbon graph as a graph
embedded in $S^3$ together with a ribbon (a long, thin rectangle)
assigned to each edge and a disk assigned to each vertex with the
incident ribbons attached at the boundary in such a way that an
oriented surface results.  A graph projection is then the projection of
the graph in the usual sense, where the normal vector everywhere points
upward out of the plane of the projection, or the surface consisting of
the ribbons and disks lies flat against the plane of the projection.  A
tangled ribbon graph with boundary can be thought of as a graph, with
its ribbons or normal vectors, embedded in a ball in $S^3$ and with its
endpoints lying at the boundary of the ball.

We define a freeway invariant to be a function defined on freeway
projections which depends only on the underlying freeway.

\subsection{Linearly recurrent invariants}

The topological invariants presented in this paper are best understood
as linearly recurrent invariants.

Let $I(G)$ be some $k$-valued invariant of some tangled graphs $G$
which are possibly decorated in some way.  There is a natural way to
extend the invariant $I$ to an invariant $\vec{I}(B)$ for graph with
boundary $B$ which takes values in a certain vector space over $k$
which depends only on $\bd B$.  If $b$ is a boundary, then we let
$F(b)$ be the vector space of all $k$-valued functions on the set of
all links with boundary $-b$.  Then if $\bd B = b$, we can choose
$\vec{I}(B) \in F(b)$ to simply be the function which takes $B'$ to
$I(B \cup B')$.  We call the resulting invariant $\vec{I}$ the
vector invariant associated to the scalar invariant $I$.

So far, the notion of $\vec{I}$ is not very interesting because the
vector space $F(b)$ may be very large.  However, it may happen that
only a small portion of $F(b)$ is needed to understand the original
invariant $I$.  We define the spanning space $\Span(b)$ to be the
span of $\vec{I}(B)$ for all $B$ such that $\bd B = b$.  We say that
$I$ is a linearly recurrent invariant, or a recurrent invariant
for short, if $\Span(b)$ is finite-dimensional for all $b$.  The
motivation for the name is that if the dimension of the space
$\Span(b)$ is small, then it is often possible to define the original
invariant $I$ recursively by describing the invariant $\vec{I}(b)$ and
the vector space it lives in for only a few choices of $B$, a strategy
which is similar to the definition of a linearly recurrent sequence of
numbers.

Another more suggestive notation for scalar and vector invariants is
bra-ket notation, by analogy with quantum mechanics.  In this notation
the value of the scalar invariant for some link or tangled graph $G$ is
$\langle G \rangle$, while the corresponding vector invariant for some
link or tangled graph with boundary $B$ is denoted $\mid B\rangle$.

As a warm-up for the rest of the paper, we review the derivation and
definition of the Kauffman bracket, a linearly recurrent invariant of
regular isotopy of links.  We denote the Kauffman bracket by
$\langle\cdot\rangle_{A_1}$ and by the bra-ket convention we denote the
corresponding vector invariant by $\mid\cdot\rangle_{A_1}$.  We assume
that the 4-point spanning space is 2-dimensional and that the two
crossingless tangles with 4 endpoints are a basis, and we assume that
the 0-point spanning space is 1-dimensional.  In this case there are
necessarily relations: $$\includegraphics{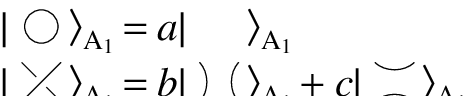}$$
If we assume that the bracket is invariant under the first of the
regular isotopy moves, we can compute: $$\includegraphics{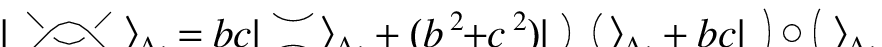}$$
to obtain the consistency equations:
\begin{eqnarray*}
bc & = & 1 \\
abc + b^2 + c^2  & = & 0
\end{eqnarray*}
We must also check the second regular isotopy move.  However, by a
ubiquitous trick in the theory of linearly recurrent invariants, if
the ket of a crossing is a linear combination of kets of crossingless
tangles (or in other contexts crossingless graphs), the second regular
isotopy move is a corollary of the first regular isotopy move (more
generally, a corollary of other moves involving vertices).  We
parameterize the solutions by setting $b = -q^{1/4}$, which is slightly
peculiar but is natural with hindsight.  We obtain the recursive rules
described in the introduction.  The invariant is completely determined
except for its value at the empty link, which lies at the base of the
recursive rules.  The choice of this value is equivalent to the choice
of a global normalization factor.  The most natural normalization in
the context of this paper is to choose $\langle\emptyset\rangle_{A_1} =
1$.

\subsection{Constructing a freeway invariant \label{secconstruct}}

We assume that there is a scalar freeway invariant
$\langle\cdot\rangle_{G_2}$ whose corresponding vector invariant is
$\mid\cdot\rangle_{G_2}$ and which satisfies the following condition:
For $n<6$, the bracket of the acyclic, crossingless freeways with $n$
endpoints are linearly independent vectors, and that the bracket of a
freeway which consists of a face with $n$ sides is a linear combination
of these acyclic freeways.  These freeways are precisely the ones
listed in the statement of Theorem~\ref{existg2}.  Thus, we assume the
existence of coefficients $a$, $b$, $c$, $d_1$, $d_2$, $e_1$, and $e_2$
as given in the statement of Theorem~\ref{existg2} (but we do not yet
assume the formulas given for those coefficients).

We call this assumption the 1,0,1,1,4,10 Ansatz because there are
1,0,1,1,4, and 10 acyclic freeways with 0,1,2,3,4, and 5 endpoints,
respectively.  Note that by symmetry there can only be two distinct
coefficients for the relation for a square and pentagon, for otherwise
the bracket of a square or a pentagon would satisfy more than one
linear relation, which would violate the linear independence
assumption.

\begin{theorem} \label{uniqg2} Every freeway invariant which satisfies
the 1,0,1,1,4,10 Ansatz is equal to the invariant defined in
Theorem~\ref{existg2} for some value of $q$.  \end{theorem}

We will proceed by proving Theorem~\ref{uniqg2}.  Along the way all of
the conditions for the existence of the invariant will be met, thereby
proving Theorem~\ref{existg2}.

The following lemma shows that the coefficients $a$,$\ldots$,$e_2$
completely determine $\langle\cdot\rangle_{G_2}$ for crossingless
freeways, provided that the coefficients are chosen so that the bracket
exists, up to a choice of global normalization.  As in the case of
the Kauffman bracket, we choose for the sake of convention
$\langle\emptyset\rangle_{G_2} = 1$.

\begin{lemma} Let $Z$ be a non-empty crossingless freeway.  The
graph $Z$ has at least one simply-connected face with five or fewer sides
which does not share an edge with itself.
\end{lemma}
\begin{proof} If we allot each face one third of each of its vertices
and one half of each of its edges, then the Euler characteristic of a
simply-connected face with $n$ sides is $1-n/6$ and the Euler
characteristic of a multiply-connected face is negative or zero.  If no
face shares an edge with itself, then since the Euler characteristic of
the sphere is 2, we must have at least one face with positive Euler
characteristic.

Suppose instead that at least one face does share an edge with itself.
Then we may draw a circle which is inside this face except at one point:
$$\includegraphics{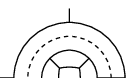}$$
If we draw a circle for each such face, then at least one of the circles
must be innermost, i.e. it must bound a region containing no such faces.
The total Euler characteristic of the faces in this region is at least
5/6, so once again there must be a face in the region with the desired
property.
\end{proof}

A variation of the proof of this lemma yields a self-justification for
the 1,0,1,1,4,10 Ansatz:  If a crossingless freeway has $n$ endpoints,
then either the freeway is acyclic or the total Euler characteristic of
its faces is at least $1-n/6$.  If $n<6$, then any freeway with $n$
endpoints is a linear combination of acyclic freeways.  The Ansatz
implies that the acyclic freeways are bases for the first six spanning
spaces of the bracket.

We now derive the values for the coefficients $a$, $b$, $c$, $d_i$, and
$e_i$ which are consistent for the values of the bracket on all
crossingless freeways.  Let $P$ be a minimal crossingless freeway
projection for which the recursive definition contradicts itself.  Suppose
that the projection $P$ has two faces $A$ and $B$ with five or fewer
sides which are not adjacent.  Then reducing the face $A$ must yield
the same result as reducing the face $B$, because after we reduce
$A$ we are still free to reduce $B$, and vice versa.  Thus, the only
possibility of a contradiction arises when $P$ has two adjacent faces
with at most five sides.

Suppose, for example, that $P$ has a triangle next to a pentagon; it
suffices to consider the case when $P$ is a triangle next to a pentagon
and nothing more, like so: $$\includegraphics{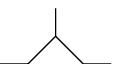}$$
We can in this case reduce either the pentagon or the triangle:
$$\includegraphics{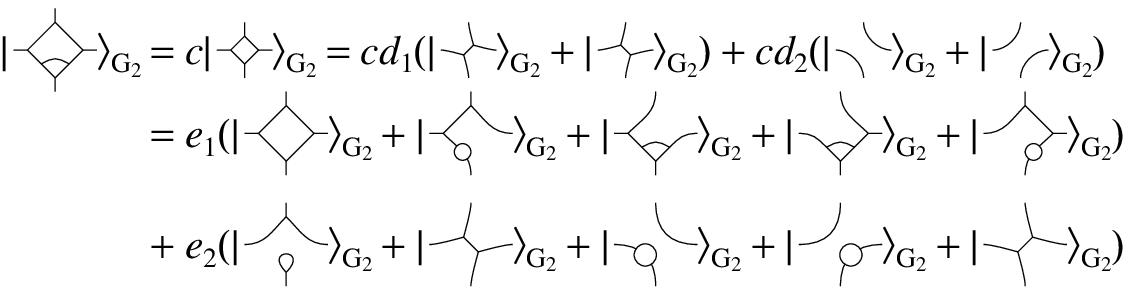}$$
By the linear independence assumptions, we obtain the scalar
equations:

\begin{eqnarray*}
c d_1 & = & e_1 d_1 + e_1 c + e_1 b + e_2 \\
c d_2 & = & e_1 d_2 + e_2 b
\end{eqnarray*}

We can perform this procedure for all cases of a face with
$n$ sides adjacent to a face with $m$ sides except $n = m = 5$.
We obtain the equations:

\begin{eqnarray}
b^2               & = & b d_1 + d_2 + a d_2 \label{fifthtodo} \\
c^2               & = & b d_1 + c d_1 + d_2 \nonumber \\
b c               & = & 2 e_1 b + 2e_2 + a e_2 + e_1 c \nonumber \\
c d_1             & = & d_1 e_1 + c e_1 + b e_1 + e_2 \label{fourthtodo} \\
c d_2             & = & d_2 e_1 + b e_2 \nonumber  \\
d_1 e_1 + d_2     & = & e_1^2 \label{secondtodo} \\
d_1 e_1 + d_1^2   & = & e_1^2 + c e_1 + d_1 e_1 \label{thirdtodo} \\
d_1 e_1           & = & e_1^2 + d_1 e_1 + e_2 \label{firsttodo} \\
d_1 e_2 + d_2 c   & = & e_1 e_2 + b e_2 + 2 d_2 e_1 \nonumber \\
d_1 e_2           & = & e_1 e_2 + d_2 e_1 \nonumber \\
d_1 e_2 + d_1 d_2 & = & e_1 e_2 + e_2 c \nonumber
\end{eqnarray}

Since we have not made any linear independence assumptions for
freeways with six endpoints, we do not obtain any equations which
the coefficients must satisfy {\em a priori} in the case of two
neighboring pentagons.  Nevertheless, we must check the equations
anyway, since the relevant freeways might be linearly independent.
Here are the equations which result:
\begin{eqnarray*}
e_2 d_1 & = & e_1 d_2 + e_1 e_2 \\
e_2     & = & e_1^2
\end{eqnarray*}

The only remaining possibility is that there are two faces with
five or fewer sides which have two or more edges as their common
border.  This case is inconsequential and the reason is left as an
exercise to the reader.

Undeterred by the sprawl of these equations and by the fact
that there are more constraints than unknowns, we can proceed by
reducing them to a sequence of relatively simple relations.  First, we note
that one solution can be varied in an uninteresting way to obtain
another:  We can multiply $b$, $c$, $d_1$, and $e_1$ by a constant and
$d_2$ and $e_2$ by the square of that constant.  To eliminate this
degree of freedom, we normalize by declaring that $e_1 = 1$.
(The case $e_1 = 0$ produces only trivial solutions.)  Then we can
eliminate variables by using some of the equations to simplify the
algebra:  We find that  $e_2 = -1$ by equation (\ref{firsttodo}),
then express $d_2$ in terms of $d_1$ by equation (\ref{secondtodo}),
then $c$ in terms of $d_1$ by equation (\ref{thirdtodo}), then
we find $b$ by equation (\ref{fourthtodo}), and finally $a$ by equation
(\ref{fifthtodo}).  In the process we see that the other equations
are satisfied automatically, and we are left with a second free parameter.
We make the useful reparameterization $d_1 = - q - q^{-1}$, with
$q$ an indeterminate, to obtain the solutions given in the statement
of Theorem~\ref{existg2}.

To extend the invariant to links and to freeways with crossings, we
further assume that a crossing is a linear combination of the four
acyclic, crossingless freeways with four endpoints.  To find the four
coefficients $f_1$, $f_2$, $g_1$, and $g_2$ for this linear dependence,
we must investigate the behavior of the bracket under the regular isotopy
moves.  We first reduce the ``before'' and ``after'' pictures of the
first regular move:
$$\includegraphics{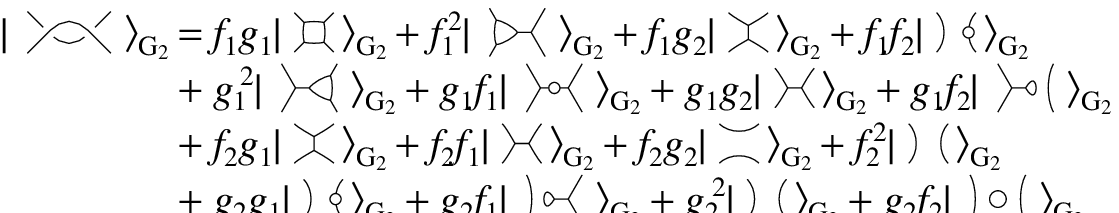}$$
This yields the equations:

\begin{eqnarray}
f_1 g_1 d_1 + f_1^2 c + g_1^2 c + f_1 g_1 b + g_1 g_2 + f_1 f_2 & = & 0
\nonumber \\
f_1 g_1 d_1 + f_1 g_2 + f_2 g_1 & = & 0 \nonumber \\
f_1 g_1 d_2 + f_1 f_2 b + f_2^2 + g_1 g_2 b + g_2^2 + f_2 g_2 a & = & 0
\nonumber \\
f_1 g_1 d_2 + f_2 g_2 & = & 1  \label{usemea}
\end{eqnarray}
the third and fourth moves result in the equations:
\begin{eqnarray}
f_1^2 d_1 + f_1 g_1 c + g_1 f_1 e_1 g_1^2 d_1 & = & 0 \nonumber \\
g_1 f_1 e_1 + g_1 f_2 & = & f_1 \nonumber \\
g_1 f_1 e_1 + g_2 f_1 & = & g_1 \nonumber \\
g_1 f_1 e_1 + g_1^2 d_1 + g_2 g_1 & = & 0 \nonumber \\
f_1^2 d_1 + f_1 f_2 + g_1 f_1 e_1 & = & 0 \nonumber \\
g_1 f_1 e_2 + g_2 f_2 & = & 0 \label{usemeb} \\
f_1^2 d_2 + f_1 g_2 b + g_1 f_1 e_2 + g_1 g_2 c + g_2^2 & = & 0 \nonumber \\
g_1 f_1 e_2 + g_1^2 d_2 + f_2 g_1 b + f_2 f_1 c + f_2^2 & = & 0 \nonumber \\
f_1^2 d_2 + g_1 f_1 e_2 & = & f_2 \label{usemec} \\
g_1 f_1 e_2 + g_1^2 d_2 & = & g_2 \label{usemed}
\end{eqnarray}

Comparing equations (\ref{usemea}) and (\ref{usemeb}), we obtain:
\bq f_1 g_1 = f_2 g_2 = \frac{1}{1+ d_2} = \frac{1}{q + 2 + q^-1}
\label{usemee} \eq
Multiplying equations (\ref{usemec}) and (\ref{usemed}) and substituting
the formula (\ref{usemee}), we obtain:
$$ f_1^2 + g_1^2 = \frac{q + q^{-1}}{q + 2 + q^{-1}} $$
Since we now know the sum and the product of $f_1^2$ and $g_1^2$, we
obtain a quadratic equation.  The result is the solutions given in
the statement of Theorem~\ref{existg2}.

It remains to check the second regular isotopy move.  As in the case of
two adjacent pentagons, we do not obtain any equations which the
coefficients must satisfy {\em a priori}, but we must check the move
anyway.  In this case the ubiquitous trick which was left as an
exercise to the reader in the case of the Kauffman bracket appears
again.  If we have invariance under the other moves, and if a crossing
is a linear combination of crossingless freeways, then the second
regular isotopy move is automatically satisfied.

In solving the quadratic equation we must break a Galois symmetry which
maps the parameter $q$ to $q^{-1}$.  A rationale for making the choice
presented here, as well as for choosing the parameter $q$ in the first
place, will be presented in the second part of the paper.

\subsection{Two other polynomial invariants}

We can repeat the approach of Section~\ref{secconstruct} two more
times by changing the notion of a freeway and picking a new Ansatz.
To avoid confusion, let us define a $G_2$ freeway to be what we 
previously called a freeway.

The solutions will be polynomial invariants of a single variable.
However, for the purpose of a uniform treatment of all of the
invariants, it is desirable in general to make this variable a certain
root of a variable $q$ which will be reused for each of the
invariants.  Therefore the final answers will be elements of the field
$\C(q,q^{1/2},q^{1/3},q^{1/4},\ldots)$.  For every rational number $a$
we distinguish a power $q^a$ of the indeterminate $q$ in such a way
that for every integer $n$, $q^{an} = (q^{a})^n$.

We define an $A_2$ freeway to be an oriented planar graph with
crossings and junctions that look like this: $$\includegraphics{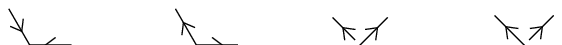}$$
We assume the existence of a regular isotopy invariant
$\langle\cdot\rangle_{A_2}$ which satisfies the 1,0,1,1,2 Ansatz:  The
acyclic, crossingless $A_2$ freeways with 0,1,2,3, or 4 endpoints are
assumed to be bases for their respective spanning spaces.  This Ansatz,
although simpler, comes with some technicalities.  In the 2-endpoint
case, one endpoint must be oriented inward and the other outward.  In
the 4-endpoint case, two endpoints must be oriented inward and two
outward, and we get two different bases which differ in the
ordering of the endpoints in the projection: $$\includegraphics{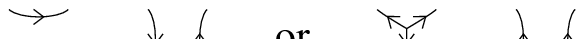}$$
The two possible boundaries are equivalent, but not canonically, so
there is no canonical formula for a change of basis.  The endpoints of
a 3-endpoint freeway must all be oriented inward or all outward, but
in fact the 3-endpoint case of the Ansatz is not needed in the
derivation of the invariant.  The 1-endpoint case is even less relevant
to the derivation because there are no $A_2$ freeways whatsoever with
only 1 endpoint.

By a derivation similar to that of Section~\ref{secconstruct}, the only
possibilities for $\langle\cdot\rangle_{A_2}$ are the ones described by
Theorem~\ref{exista2}.  As before, the rules given in the theorem
suffice as an recursive definition because any trivalent planar graph
must have a face with at most 5 sides, except the situation is simpler
because all faces of an $A_2$ freeway must have an even number of
sides.  The invariant $\langle\cdot\rangle_{A_2}$ is identically equal
to $RT_{A_2}$ if all edges are colored with the 3-dimensional
representation $V_{1,0}$.  The edges are oriented because the
representation is not self-dual.  The invariant is a special case of
the HOMFLY polynomial with a normalization that makes it a regular
isotopy invariant rather than an isotopy invariant:  $$\includegraphics{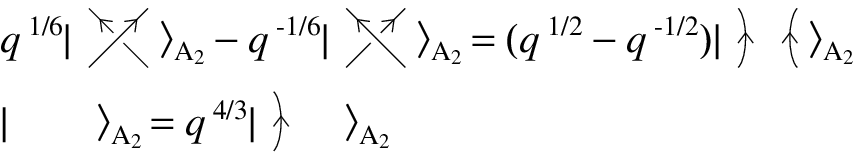}$$

We also define a $C_2$ freeway to be a planar graph with two
kinds of edges, illustrated as double and single edges, with junctions
and crossings that look like this: $$\includegraphics{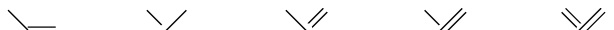}$$
The Ansatz for $\langle\cdot\rangle_{C_2}$ has a somewhat different form.
We first assume the existence of a relation which we
call a switching relation: $$\includegraphics{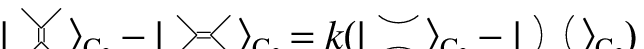}$$

Associated to the switching relation is a switching move: $$\includegraphics{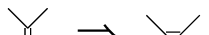}$$
We define a switching class of $C_2$ freeways to be an equivalence class
under switching moves.  We say that a double edge of a $C_2$ freeway
is external if at least one of its vertices is an endpoint of the
freeway, otherwise the edge is internal.  Then the Ansatz says that when
$s+2d < 8$, a collection of acyclic, extended $C_2$ freeways with $s$
single endpoints and $d$ double endpoints and no internal double edges
are a basis if the collection has one representative from each
switching class.  The caveat about ordering of endpoints applies here
just as in the $A_2$ case.

By a third iteration of Section~\ref{secconstruct}, the only possibilities
for $\langle\cdot\rangle_{C_2}$ are those of Theorem~\ref{existc2}.
To show that these rules suffice as an inductive definition,
we define the angles of a (trivalent) junction to be $3\pi/4$, $3\pi/4$,
and $\pi/2$, like so: $$\includegraphics{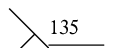}$$
We define the curvature of a face to be $2\pi$ minus the sum of the
exterior angles.  By an Euler characteristic argument, every
crossingless $C_2$ freeway without boundary must have a face with
positive curvature, and such a face can be reduced by the switching
move to one for which one of the other moves applies.  

Another approach to the $C_2$ bracket is to define an extended
$C_2$ freeway to be one which can also have tetravalent vertices with
no incident double edges.  A tetravalent freeway is one with no
double edges whatsoever.  We assume the existence of a relation:
$$\includegraphics{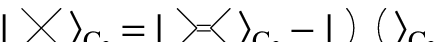}$$
where the constants $\alpha$ and $\beta$ are arbitrary except that
$\alpha$ must be non-zero to avoid degeneracy.  If we interpret the
four angles at a tetravalent vertex to be 90 degrees, we can derive
recurrence rules for tetravalent freeways by reducing faces with
positive curvature, namely faces with at most 3 sides.  A face with $n$
sides now has $2n$ endpoints, and for $n = 0,1,2,3$ the number of
acyclic tetravalent graphs with $2n$ endpoint is 1,1,3, and 14.  But by
the Ansatz for the $C_2$ bracket, these acyclic graphs must be a basis
provided that $\alpha \ne 0$.  Thus we obtain an example of a linearly
recurrent theory for tetravalent graphs; double edges and trivalent
vertices are avoided.  (Note that a tetravalent vertex should not be
confused with a crossing.)  Moreover, some of the choices of the
constants $\alpha$ and $\beta$ result in a simpler recurrent calculus
than others.  A particularly fortuitous choice is $\alpha = 1/2$ and
$\beta = 1/2$, whose first advantage is a simpler relationship between the
trivalent and tetravalent theories: $$\includegraphics{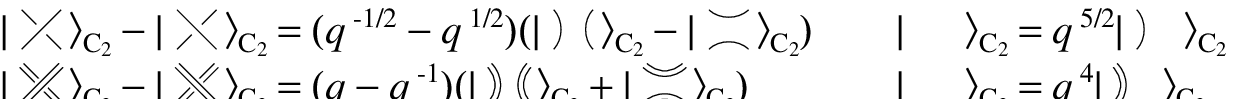}$$
Its second advantage is that in the recurrence relations, the
coefficients of 1 of the acyclic graphs with 4 endpoints and 10 of the
acyclic graphs with 6 endpoints vanish.  The resulting recurrence
relations are those given in the second half of the statement of
Theorem~\ref{existc2}.

\begin{question} Currently the tetravalent $C_2$ bracket, for various
values of $\alpha$ and $\beta$, merely represents examples of linearly
recurrent invariants that satisfy the 1,1,3,14 Ansatz.  What are all
the invariants of tangled tetravalent graphs that satisfy this Ansatz?
\end{question}

The $C_2$ bracket constitutes two special cases of the Kauffman
polynomial, one for plain edges and one for double edges.  The two
relations it satisfies, both in ``Dubrovnik'' form, are:
$$\mbox{(Figure 27)}$$

\subsection{A few calculations}

In this section we present the results of a few calculations with
$\langle\cdot\rangle_{G_2}$.  We begin with the auxiliary recursive
rules:
$$\mbox{(Figure 28)}$$
Note that these rules can be inferred from
\cite[p. 61]{Reshetikhin:g2}, and that they allow one to compute
the bracket for any knot with 7 or fewer crossings.  The simplest link
which seems to defeat these rules is the Borromean rings, with 6
crossings and 3 components.  Unfortunately, the Borromean rings have
also defeated the author's attempts at calculations by hand.  However,
some other possibly interesting values of the bracket and relations of
the ket follow:
$$\mbox{(Figure 29)}$$
Note that all of the links and graphs are evaluated with the projection
framing.

\section{Invariants of Lie algebra representations \label{seccompare}}

In this part of the paper we place the invariants
$\langle\cdot\rangle_{\mathfrak{g}}$ in the proper context of the invariants
$RT_{\mathfrak{g}}$ defined by Reshetikhin and Turaev in
\cite{RT:ribbon} and study their properties from the point of
view of representation theory of Lie algebras.

\subsection{Invariants of colored ribbon graphs}

We begin with a summary of Reshetikhin and Turaev's definition and
some basic properties.

We review some of the theory of finite-dimensional left modules, or
representations, of a Hopf algebra $H$.  If we have two representations
$V$ and $W$, the tensor product $V \tensor W$ is also a representation
by the coalgebra structure of $H$.  There is a trivial representation
$I$ such that $I \tensor V \cong V \tensor I \cong V$ canonically.
However, there is no canonical isomorphism between $V \tensor W$ and $W
\tensor V$ as representation of $H$, and indeed the two modules need
not be isomorphic at all.  In particular, the switching map, which
takes $x \tensor y$ to $y \tensor x$, is in general not a module
homomorphism.  On the other hand, we do have $U \tensor (V \tensor W)
\cong (U \tensor V) \tensor W$, and it is still true that if $A$ and
$B$ are morphisms between representations, $A \tensor B$ is also.
Also, for each representation $V$ there is a representation $V\dual$
and a canonical identification of $\Hom(U \tensor V\dual,W)$ and
$\Hom(U,W \tensor V)$.  There is also a representation $V\odual$ with a
canonical identification of $\Hom(V\odual \tensor U,W)$ with $\Hom(U,V
\tensor W)$.  Although there are canonical isomorphisms $V\odual\dual
\cong V\dual\odual \cong V$, $V\dual\dual$ need not be isomorphic to
$V$.  Finally, we define $\Inv(V) = \Hom(I,V)$ by analogy with the case
when $H$ is a group algebra.

There is a visual notation for morphisms in the category of $H$-modules
and their compositions.  If $f$ is a morphism from $V_1 \tensor \ldots
\tensor V_n$ to $W_1 \tensor \ldots W_n$, we can represent $f$ by
a diagram with one vertex:
$$\mbox{(Figure 30)}$$
Similarly for compositions of morphisms.  For example, if
$f:V_1 \to V_3 \tensor V_4$, $g:V_4 \tensor V_2 \to V_5$,
and $h:V_3 \tensor V_5 \to V_6 \tensor V_7$ are morphisms, then
we can express $h \circ (I \tensor g) \circ (f \tensor I)$
by the diagram: $$\mbox{(Figure 31)}$$
Note that the facts stated in the previous paragraph allow
us to erase internal edges labeled with the representation $I$; by
a slight abuse of notation we may, for example, draw a morphism
$f:I \to V \tensor W$ as: $$\mbox{(Figure 32)}$$

In general these diagrams must be planar graphs in which none of the edges pass
through a slope of zero.  However, it is sometimes possible to augment the
structure of a Hopf algebra to obtain a ribbon Hopf algebra, in which the
graphs can have edges that bend arbitrarily and can have edges that cross each
other, provided one arm of each crossing is labeled ``over'' and the other
``under''.  First, there is a distinguished element $R$ of $H \tensor H$ which
we can use to construct, for each $V$ and $W$, an isomorphism $\breve{R}:V
\tensor W \to W \tensor V$ which is denoted by a right-handed crossing:
$$\mbox{(Figure 33)}$$
The inverse $\breve{R}^{-1}$, which need not equal $\breve{R}$, is denoted by a
left-handed crossing: $$\mbox{(Figure 34)}$$ Second, there is a distinguished
element $K$ of $H$ which we can use to construct an isomorphism
$\breve{K}:V\dual\dual\to V$.  It follows that there exist canonical morphisms
which can be denoted as follows:
$$\mbox{(Figure 34)}$$
In light of the existence of these morphisms and the identification of $V\dual$
and $V\odual$, we can orient the edges of the diagrams and declare that an edge
labeled with $V\dual$ is equivalent to an oppositely oriented edge labeled with
$V$.  In a ribbon Hopf algebra the elements $R$ and $K$ must satisfy certain
axioms which can be translated to various identities for the corresponding
morphisms. Here are some of these identities:
$$\mbox{(Figure 35)}$$

These identities allow us to interpret the finite-dimensional representation
theory of a given ribbon Hopf algebra $H$ as an invariant of colored ribbon
graphs, where each edge is colored by a representation of $H$ and each vertex
is colored by an element of $\Inv(V_1 \tensor V_2 \tensor \ldots V_n)$, where
$V_1,\ldots,V_n$ are the colors of the incoming edges going clockwise around
the vertex. (The vertex must also be marked to identify which edge corresponds
to the first factor in the tensor product; however there is an identification
between the sets of colors allowed for two different markings.  Although there
is no single preferred isomorphism between $V \tensor W$ and $W \tensor V$,
there is such an isomorphism between $\Inv(V \tensor W)$ and $\Inv(W \tensor
V)$.  It is constructed from the $\breve{K}$ morphisms.) We may interpret the
value of a ribbon graph with no endpoints as a morphism from $I$ to $I$, i.e.
as a scalar invariant.  But the invariant is also a linearly recurrent
invariant; the spanning space of endpoints colored with $V_1,\ldots,V_n$ is
again $\Inv(V_1 \tensor V_2 \tensor \ldots V_n)$.  This invariant of colored
ribbon graphs is the Reshetikhin-Turaev invariant corresponding to the ribbon
Hopf algebra $H$.

We remark that for a fixed graph, the invariant is multilinear under
direct sums, and more generally under extensions, of the colors of the
edges, so it suffices to consider irreducible representations as
colors.  We also consider the consequences of an irreducible
representation $V$ being self-dual.  If $V$ is symmetrically self-dual,
i.e. there is an element $B$ of $\Inv(V \tensor V)$ which is invariant
under the switching map, then we can allow an unoriented edge colored
by $V$ by interpreting it as two oriented edges using $B$:
$$\mbox{(Figure 36)}$$
On the other hand, if $V$ is anti-symmetrically self-dual, i.e.
the invariant tensor $B$ negates under the switching map, then this
construction has a sign ambiguity.  But suppose that all representations
of $H$ are self-dual, some of them anti-symmetrically self-dual.
Then we may resolve the ambiguity by observing that the Reshetikhin-Turaev
theory works equally well for Hopf superalgebras.  If we reinterpret
anti-symmetrically self-dual representations as being negatively graded,
they become supersymmetric.  It is possible to replace the Hopf algebra
$H$ with a Hopf superalgebra $H'$ for which these are the natural
gradings.

The other essential ingredient in the Reshetikhin-Turaev theory is a
collection of examples of ribbon Hopf algebras.  For each simple Lie
algebra ${\mathfrak{g}}$ there is a Hopf algebra $U_q({\mathfrak{g}})$, the
quantized universal enveloping algebra of ${\mathfrak{g}}$, over the field
$\C(q,q^{1/2},q^{1/3},\ldots)$ (see \cite{Drinfeld:quantum}).  The representation
ring of this Hopf algebra is the same as that of ${\mathfrak{g}}$ or
$U({\mathfrak{g}})$.  Consequently we will refer to the irreducible
representations $U_q({\mathfrak{g}})$ as $V_\lambda$, where $\lambda$ is the
highest weight of the corresponding representation of ${\mathfrak{g}}$.
However, the calculus of morphisms between the representations is only
the same as that for ${\mathfrak{g}}$ at the specialization $q=1$.  We
denote the corresponding ribbon graph invariant by $RT_{\mathfrak{g}}$.

The most important first case of the ribbon graph invariants is their
value on the untwisted unknot colored with some $V_\lambda$.  The value
of $RT_{\mathfrak{g}}$ for this graph is known as the quantum dimension of
$V_\lambda$ and equals the trace of the element $K$.  This element is a
member of the subalgebra of $U_q({\mathfrak{g}})$ generated by the Cartan
subalgebra of ${\mathfrak{g}}$, and this subalgebra is the same as the
corresponding subalgebra of a certain completion $U({\mathfrak{g}})$, the
ordinary universal enveloping algebra.  If $\rho$ is half of the sum of
the positive roots of the root system of ${\mathfrak{g}}$ and $H_\rho$ is the
corresponding element of the Cartan subalgebra, then $K = q^{H_\rho}$.
The upshot is that the trace of $K$ on $V_\lambda$ is $\sum_\alpha
m(\alpha) q^{(\rho,\alpha)}$, where $\alpha$ runs over the weights
of $V_\lambda$ and $m(\alpha)$ is the multiplicity of $\alpha$.
Another formula for the quantum dimension is the quantum Weyl
dimension formula:
$$\dim_q V_\lambda = \frac{\prod_{\alpha \succ 0} \left[(\lambda+\rho,\alpha)
\right]}{\prod_{\alpha \succ 0}
\left[(\rho,\alpha) \right]}, $$
where $[n] = \frac{q^{n/2}-q^{-n/2}}{q^{1/2}-q^{-1/2}}$ is a ``quantum
integer'' (see \cite[ch. 6, p. 151]{Bourbaki:lie} for a proof that this is
the trace of $K$).  For example, when ${\mathfrak{g}} = G_2$, we obtain:
$$\dim_q V_{1,0} = q^5+q^4+q+1+q^{-1}+q^{-4}+q^{-5} =
\frac{[7][12][2]}{[6][4]}$$
as the quantum dimension of the smallest non-trivial irreducible
representation.  Note, however, that when we perform a supersymmetric
conversion of a ribbon Hopf algebra $H$ as described above, the
quantum dimension of a negatively graded representation negates.

The next most important computation is the evaluation of a right-handed
twist: $$\mbox{(Figure 37)}$$
where $c(V_\lambda)$ is the Casimir number for $V_\lambda$.
A derivation of this equation is given in \cite{Reshetikhin:g2}.

\subsection{Comparison with $\langle\cdot\rangle_{\mathfrak{g}}$}

We restrict our attention to the exceptional simple Lie algebra $G_2$.
All representations of $G_2$ are symmetrically self-dual.  The vector
space $\Inv(V_{1,0}^{\tensor 3})$ is 1-dimensional; let us choose some
element $T$ in it.  We interpret a $G_2$ freeway as a colored
ribbon graph by coloring all edges with the 7-dimensional
representation $V_{1,0}$ and all vertices with $T$.  In this way we
obtain a value of $RT_{G_2}$ for $G_2$ freeways.

We can verify by a computation with weight diagrams that $\dim
\Inv(V_{1,0}^{\tensor n})$ is 1,0,1,1,4, and 10 as $n$ ranges from 0 to
5 (see \cite{Humphreys:gtm}).  A lengthier calculation shows that the values
of $RT_{G_2}$ on acyclic graphs with $n$ endpoints are indeed a basis
for $n<5$ and $q=1$ (and therefore also for $q$ an indeterminate).  The
author and Rena Zieve performed this calculation explicitly.  The
tensor $T$ was chosen as a trilinear form $T(a,b,c)$ on $V_{1,0}$.  It
is well-known that for the compact real Lie group $G_2$, $V_{1,0}$ can
be interpreted as the imaginary octonions, or Cayley numbers, which are
equipped with a positive-definite inner product and a cross product, a
bilinear function ``$\times$'' with the property that if $a$ and $b$ are
orthogonal unit vectors, $a \times b$ is a unit vector which is
orthogonal to both $a$ and $b$.  We define $T(a,b,c) = (a \times b)
\cdot c$.

The conclusion is that the invariant $RT_{G_2}$ satisfies the
1,0,1,1,4,10 Ansatz.  It remains only to show that the parameter $q$ in
$RT_{G_2}$ agrees with the parameter $q$ in
$\langle\cdot\rangle_{G_2}$.  To do this we first check
that the quantum dimension of $V_{1,0}$ equals the constant $a$
in the recurrence relation of $\langle\cdot\rangle_{G_2}$, and we
compute the factor gained by a right-hand twist: $$\mbox{(Figure 38)}$$
We verify that $c(V_{1,0}) = 6$.  Thus, the two parameterizations
are the same, up to a constant factor in the choice of $T$.
(The choice of $T$ above for the explicit $q=1$ computation
agrees with the right choice for $\langle\cdot\rangle_{G_2}$.)
Notice in particular that the Casimir element is always positive,
which assures us that the ambiguity between $q$ and $q^{-1}$ has
been resolved in agreement with convention.

By similar but simpler reasoning we can show that
$\langle\cdot\rangle_{A_2} = RT_{A_2}$, where all edges of an $A_2$
freeway are colored with the 3-dimensional representation $V_{1,0}$,
whose dual is $V_{0,1}$.  The orientations of the edges are needed
because $V_{1,0}$ is not self-dual.  The colors for the vertices can be
recognized as the determinant, or the usual 3-dimensional cross
product.

We also have $\langle\cdot\rangle_{C_2} = RT_{C_2}$, where the single
edges are colored with the 4-dimensional representation $V_{1,0}$ and
the double edges are colored with the 5-dimensional representation
$V_{0,1}$.  The representation $V_{1,0}$ is the canonical
representation of $\sp(4) \cong C_2$.  Note that it is
anti-symmetrically self-dual.  To convert from
$\langle\cdot\rangle_{C_2}$, where edges are not oriented, we must use
the Hopf superalgebra $U_q(\sp(4))'$ instead of $U_q(\sp(4))$.  This
superalgebra can also be thought of as $U_q(\osp(0,4))$, where
$osp(0,4)$ is an orthosymplectic Lie superalgebra.  (Similarly, the
Kauffman bracket $\langle\cdot\rangle_{A_1}$ more properly describes
invariants from the representation theory of $\osp(0,2)$ rather than the
theory of $\sl(2) = \sp(2)$.)  Finally, the representation $V_{0,1}$ is
the canonical representation of $\so(5) \cong B_2 \cong C_2$.

\subsection{Questions open to the author}

The coefficients $a$, $b$, $c$, $d_1$, $d_2$, $e_1$, and $e_2$ in the
definition of the $G_2$ bracket, as well as their analogues for the
$A_2$ and $C_2$ brackets, are all Laurent polynomials in $q$ with
either nonpositive or nonnegative integer coefficients.  This is a
coincidence which cannot be explained by the freedom in choosing the
parameter $q$.  The coefficients could be complicated rational
or even algebraic functions {\em a priori}.  Moreover, using the
notion of a quantum integer given above, we can write:

\begin{eqnarray*}
a & = & \frac{[12][7][2]}{[6][4]} \\
b & = & -\frac{[8][3]}{[4]} \\
c & = & \frac{[6]}{[2]} \\
d_1 & = & -\frac{[4]}{[2]} \\
d_2 & = & [3] \\
e_1 & = & 1 \\
e_2 & = & -1
\end{eqnarray*}

The coefficients for the $A_2$ and $C_2$ brackets have similar
expressions.  In general, saying that a Laurent polynomial in $q$
invariant under $q\mapsto q^{-1}$ with rational coefficients is a ratio
of products of quantum integers is equivalent to saying that the zeroes
of the polynomial are all roots of unity.

\begin{question}  Why do the coefficients in the definition of
$\langle\cdot\rangle_{A_2}$, $\langle\cdot\rangle_{C_2}$, and
$\langle\cdot\rangle_{G_2}$ have such a nice form?
\end{question}

The known properties of quantum dimension provide the answer
to this question in the case of the coefficient $a$ as well as
the value of the unknot for the other brackets.  However, the author
knows no analogous reasoning for the other coefficients.

Let us say that a crossingless $G_2$ freeway (also an
$A_2$ freeway) has non-positive curvature if all interior faces have at
least six sides.  From the definition of $\langle\cdot\rangle_{G_2}$,
it is clear that the set of freeways of non-positive curvature spans
the spanning space of tangles with $n$ endpoints.  It is not so clear
when $n>5$ that the same freeways span the corresponding spanning space
for $RT_{G_2}$, which can be identified with the vector space
$\Inv(V_{1,0}^{\tensor n})$.  Nevertheless, it follows from the
fundamental theorem of invariant theory for $G_2$ \cite{Schwarz:g2}:

\begin{theorem} (Schwarz) The set of all multilinear invariant
forms on the imaginary octonions is generated by addition and multiplication
of an invariant bilinear form $B(a,b)$, an invariant trilinear 
form $T(a,b,c)$ (both unique up to a scalar multiple), together
with the Hodge dual of $T(a,b,c)$.
\end{theorem}

Let us call the third form $Q$.  Then the theorem says that graphs of this
form: $$\mbox{(Figure 39)}$$
span $\Inv(V_{1,0}^{\tensor n})$.  Strictly speaking, the theorem only
applies when $q=1$, but by genericity, it must also hold for $q$ an
indeterminate.  Also, although the graphs mentioned by the theorem have
crossings and also a new tetralinear form, both of these are some
linear combination of the four acyclic $G_2$ freeways with four
endpoints.  Thus, the set of all freeways span, and therefore the ones
with non-positive curvature do too.

\begin{conjecture}  The set of $G_2$ freeways of non-positive curvature
with $n$ endpoints is a basis for $\Inv(V_{1,0}^{\tensor n})$ via
the function $\langle\cdot\rangle_{G_2}$.
\end{conjecture}

It suffices to show that the number of such freeways equals
$\dim \Inv(V_{1,0}^{\tensor n})$.  The latter number is also equal
to the coefficient of the $x^2y^3$ term in the polynomial:
\begin{eqnarray*}
(1+x+y+xy+x^{-1}+y^{-1}+(xy)^{-1})^n(x^{2}y^{3}-xy^{3}+x^{-1}y^{2}-x^{-2}y\\
+x^{-3}y^{-1}-x^{-3}y^{-2}+x^{-2}y^{-3}-x^{-1}y^{-3}+xy^{-2}-x^{2}y^{-1}+
x^{3}y-x^{3}y^{2})
\end{eqnarray*}
by the character formulas in \cite{Humphreys:gtm}.

The author has verified the conjecture
up to $n=9$ by this method; the number of acyclic, crossingless
freeways with $n$ endpoints for $n$ from 0 to 9 is 1, 0, 1, 1, 4, 10,
35, 120, 455, and 1728, respectively.

The situation is similar for $A_2$ and $C_2$.  It follows from the
classical fundamental theorem of invariance theory for $\sl(3)$ and $\sp(4)$
that freeways with non-positive curvature span.  However, the
endpoints of an $A_2$ freeway must be labeled as either ``incoming'' or
``outgoing'' depending on the orientation of the edges containing the
points, and the endpoints of a $C_2$ freeway are labeled as ``single'' or
``double''.  Recall that in the $C_2$ case we have extra equivalences 
which are conveyed by switching moves.

\begin{conjecture}  Consider a fixed set of $k$ incoming and $n$ outgoing
points on a circle.  The set of $A_2$ freeways of nonpositive curvature
having these points as endpoints is a basis for $\Inv(V_{1,0}^{\tensor n}
\tensor V_{0,1}^{\tensor k})$.  It suffices to show that the number of
such freeways equals the coefficient of the $xy^2$ term in the polynomial:
$$(xy+x^{-1}+y^{-1})^k(x^{-1}y^{-1}+x+y)^n
(xy^2-x^{-1}y+x^{-2}y^{-1}-x^{-1}y^{-2}+xy^{-1}-x^2y)$$
\end{conjecture}

\begin{conjecture}  Consider a fixed set of $k$ points labeled
``single'' and $n$ points labeled ``double'' on a circle.  A collection
of $C_2$ freeways of nonpositive curvature, with one member in each
switching class, constitutes a basis for $\Inv(V_{1,0}^{\tensor k}
\tensor V_{0,1}^{\tensor n})$.  It suffices to show that the number of
switching classes equals the coefficient of the $xy^2$ term in the
polynomial:
\begin{eqnarray*}
(x+y+x^{-1}+y^{-1})^k(1+xy+x^{-1}y+xy^{-1}+x^{-1}y^{-1})^n \\
(xy^2-x^{-1}y^2+x^{-2}y-x^{-2}y^{-1}+x^{-1}y^{-2}-xy^{-2}+x^2y^{-1}-x^2y)
\end{eqnarray*} \label{conjc2}
\end{conjecture}

Note that these conjectures have an analogue for the Kauffman bracket
as well.  In this case we may conjecture that crossingless tangles with
$2n$ endpoints form a basis for $\Inv(V_1)^{\tensor 2n}$.  By the
fundamental theorem of invariant theory they must form a spanning set;
it is well-known that the both the number of tangles and the dimension
of the invariant space equal the Catalan number $\frac{1}{n+1} \binom{2n}{n}$.

Richard Stanley and John Stembridge have shown the author a proof of
conjecture \ref{conjc2} in the case $n=0$ by a variant of the Berele
insertion algorithm \cite{Stanley:personal}, which plays the same role in the
combinatorics of $\sp(2d)$ that the Robinson-Schensted algorithm does in
the combinatorics of $\sl(d)$.  Recall that when $n=0$ we can consider
tetravalent freeways instead of switching classes of trivalent
freeways.

\begin{proof}  (Sketch) Let $F$ be a tetravalent freeway with
nonpositive curvature.  Let $v_1,v_2,\ldots,v_{2n}$ be the $2n$
endpoints of $F$.  We interpret $F$ as a collection of crossing line
segments whose endpoints are the $v_i$'s, so that the $v_i$'s are
paired by arcs of the freeway.  We let $m(i) = j$ if $v_i$ is connected
to $v_j$.  Then the function $m$ is a perfect matching, a
fixed-point-free involution of the numbers from 1 to $2n$.  Let $M(F)$
to be the matching $m$, constructed using $F$.

It is easy to show by induction on $n$ that the matching $M(F)$
satisfies what we will call the 6-point condition:  For every three
numbers $n_1$, $n_2$, and $n_3$, we cannot simultaneously have
$M(F)(n_1) > M(F)(n_2) > M(F)(n_3) > n_1 > n_2 > n_3$.  Moreover, $m$
is a bijection between freeways with non-positive curvature and
matchings that satisfy the 6-point condition.

For each $i$ from 0 to $2n$, we construct a Young tableau $T_i$, an infinite
array whose entries are positive integers or infinity.  All but
finitely-many entries of a Young tableau must be infinity, and the rows
and columns must be non-decreasing.
Here is an example of a Young tableau with the infinities omitted:
$$ \begin{array}{cccccc}
1 & 2 & 4 & 7 \\
3 & 8 \\
6 & 9
\end{array} $$
The Young tableaux $T_0$ and $T_{2n}$ are both the empty tableau, the
one with no finite entries.  We proceed to construct $T_i$ from
$T_{i-1}$ by considering $M(F)(i)$.  If $M(F)(i) > i$, we insert
$M(F)(i)$ into $T_{i-1}$ by Schensted insertion:  We replace the
highest entry $e$ in the first column which is greater than $M(F)(i)$
with $M(F)(i)$.  If $e$ is finite, we then insert $e$ by Schensted
insertion into the tableau consisting of the second column onward; if
$e$ is infinity, we stop.  If $M(F)(i) < i$, we remove $i$ by Berele
deletion:  We replace $i$ in the tableau by infinity, and then we
repeatedly switch the misplaced infinite entry by the entry below or
the entry to the right, whichever is smaller, until the array becomes a
Young tableau again.  Let $F_i$ be the support of $T_i$, the set of
place where $T_i$ is non-zero.  The sequence of $F_i$'s is called an
up-down tableau, and the algorithm described here is a bijection
between matchings of $2n$ elements and up-down tableaux of length
$2n$.  See \cite{Sundaram:cauchy} for a proof of this fact.

Matchings which happen to satisfy the 6-point condition
correspond to tableaux with at most two rows.  To see this, we observe
that if some $T_i$ has two rows and the entry in the first column
of the second row is $j$, then $j$ is the smallest number such
that $M(F)(j) < i$ and for some $k>j$, $M(F)(k) < M(F)(j)$.
Therefore if we consider the first $T_i$ with three rows, we
see that only a violation of the 6-point condition could have
resulted in transforming $T_{i-1}$ to $T_i$.  Conversely, if 
a matching does violate the 6-point condition, the corresponding
up-down tableau must at some point have three rows by the same
principle.

If an element of an up-down tableau has at most two rows, then that
element can be described by a pair of integers $(a,b)$ with $a\ge b \ge
0$, where $a$ is the number of elements in the first row and the $b$ is
the number in the second row.  Thus, we obtain, for each matching, a
sequence of pairs of integers $(a_i,b_i)$ with $a \ge b \ge 0$ such
that $a_0 = b_0 = a_{2n} = b_{2n} = 0$ and such that either $a_i$ or
$b_i$ differs by 1 from $a_{i-1}$ or $b_{i-1}$ and the other is equal.
This sequence is a lattice path in the Weyl chamber of the Lie algebra
$C_2$, and by the character theory of Lie algebras, the number of such
paths is the dimension of $\Inv(V_{1,0}^{\tensor 2n})$.
\end{proof}

We conclude with the question about the present work that most
bothers the author:

\begin{question} If ${\mathfrak{g}}$ is a simple Lie algebra of rank
greater than 2, does the invariant $RT_{\mathfrak{g}}$ have a
inductive definition based on linear recurrence, planar graphs,
and Euler characteristic?
\end{question}

\nocite{Sweedler:hopf}
% \bibliography{qa,books}

\providecommand{\bysame}{\leavevmode\hbox to3em{\hrulefill}\thinspace}

\end{document}